\documentclass{article}
\usepackage{amsmath, amssymb,latexsym}										
\begin{document}
\title{Some Relationships and Properties\\ 
       of the Hypergeometric Distribution}
\author{Peter~ H. Peskun, Department of Mathematics and Statistics \\
York University, Toronto, Ontario M3J 1P3, Canada \\ 
E-mail: peskun@pascal.math.yorku.ca}
\date{}
\maketitle
											
\begin{abstract}
The binomial and Poisson distributions have interesting relationships with the beta 
and gamma distributions, respectively, which involve their cumulative distribution 
functions and the use of conjugate priors in Bayesian statistics. We briefly 
discuss these relationships and some properties resulting from them which play an 
important role in the construction of exact nested two-sided confidence intervals 
and the computation of two-tailed P-values. The purpose of this article is to show 
that such relationships also exist between the hypergeometric distribution and 
a special case of the Polya (or beta-binomial) distribution, and to derive some 
properties of the hypergeometric distribution resulting from these relationships.

\bigskip

\end{abstract}

\begin{center}
  \textbf{1. \quad INTRODUCTION}
\end{center}

The binomial and Poisson distributions have interesting relationships with the beta 
and gamma distributions, respectively, which involve their cumulative distribution 
functions and the use of conjugate priors in Bayesian statistics. We will briefly 
discuss these relationships and some properties resulting from them in Sections 2 
and 3 for the binomial and Poisson distributions, respectively. The resulting 
properties play an important role in the construction of exact nested two-sided 
binomial and Poisson confidence intervals, and the computation of exact two-tailed 
binomial and Poisson P-values.

The purpose of this article is to show that such relationships also exist 
between the hypergeometric distribution and a special case of the Polya (or beta-binomial) 
distribution, and to derive some properties of the hypergeometric distribution 
resulting from these relationships. We shall do this in Section 4.

\begin{center}
  \textbf{2. \quad RELATIONSHIPS AND PROPERTIES OF THE BINOMIAL DISTRIBUTION}
\end{center}

Suppose that random variable $X$ has a binomial distribution with parameters $n$ 
and $p$, denoted by $X \sim \text{BIN}(n,p)$, where $n$ is a positive integer and 
$0 \leq p \leq 1$. Then, for a given $n$ and for $0 < p < 1$, the probability mass 
function (pmf) of $X$, denoted by $f_{X}(x \mid p)$, is 
\begin{align*} 
  f_{X}(x \mid p) = P(X=x \mid p) &= \binom{n}{x}p^{x}(1-p)^{n-x}, \quad x = 0,1, 
\ldots ,n, \\
                                  &=0, \quad \text{otherwise,}
\end{align*}
and $f_{X}(0 \mid 0) = f_{X}(n \mid 1) = 1$.

Suppose that random variable $Y$ has a beta distribution with parameters $\alpha > 
0$ and $\beta > 0$, denoted by $Y \sim \text{BETA}(\alpha,\beta)$. Then the 
probability density function (pdf) of $Y$, denoted by \linebreak 
$f_{Y}(y \mid \alpha,\beta)$, is 
\begin{align*}
  f_{Y}(y \mid \alpha,\beta) &= 
\frac{\Gamma(\alpha+\beta)}{\Gamma(\alpha)\Gamma(\beta)}y^{\alpha-1}(1-y)^{\beta-1}, 
\quad 0 \leq y \leq 1,\\
                             &= 0, \quad \text{otherwise,}
\end{align*}
where the gamma function $\Gamma(\kappa) = \int_{0}^{\infty} t^{\kappa-1}e^{-t} 
\, dt$ for all $\kappa > 0$.

Successive integration by parts leads to a relationship between the cumulative 
distribution functions (cdf's) of the binomial and beta distributions. If $X \sim 
\text{BIN}(n,p)$ and $Y \sim \text{BETA}(i+1,n-i)$ for integer $i$, $0 \leq i \leq 
n-1$, then 
\begin{equation} \label{E:1}
  \sum_{x=0}^{i}\binom{n}{x}p^{x}(1-p)^{n-x} = 
1-\frac{n!}{i!(n-i-1)!}\int_{0}^{p}t^{i}(1-t)^{n-i-1} \, dt.
\end{equation}
That is, $F_{X}(i \mid p) = P(X \leq i \mid p) = 1-P(Y \leq p \mid i+1,n-i) = 
1-F_{Y}(p \mid i+1,n-i)$. For fixed integer $i$, $0 \leq i \leq n-1$, it follows 
from equation (\ref{E:1}) that the function $P(X \leq i \mid p)$ is continuous and 
decreasing in $p$; for fixed integer $j$, $1 \leq j \leq n$, $P(X \geq j \mid p) = 
1-P(X \leq j-1 \mid p)$ is continuous and increasing in $p$; and for fixed integers 
$i$ and $j$, $1 \leq i \leq j \leq n-1$, $P( i \leq X \leq j \mid p)$ is 
continuous, and increasing for $0 \leq p < p_{n}(i,j)$ and decreasing for 
$p_{n}(i,j) \leq p \leq 1$ with maximum at $p = p_{n}(i,j) = \{1+[(n-i) \cdots 
(n-j)/j \cdots i]^{1/(j-i+1)}\}^{-1}$. Also, $p_{n}(0,j) = 0$ for $0 \leq j \leq 
n-1$ and $p_{n}(i,n) = 1$ for $1 \leq i \leq n$.

Suppose that the binomial parameter $p$ is unknown and we wish to estimate it. In 
Bayesian statistics, information obtained from the data x, a realization of $X \sim 
\text{BIN}(n,p)$, is combined with prior information about $p$ that is specified in 
a ``prior distribution'' with pdf $g(p)$ and summarized in a ``posterior 
distribution'' with pdf $h(p \mid x)$ which is derived from the joint distribution 
$f_{X}(x \mid p)g(p)$, and according to Bayes formula is 
\begin{equation} \label{E:2}
  h(p \mid x) = \frac{f_{X}(x \mid p)g(p)}{\int_{0}^{1}f_{X}(x \mid p)g(p) \, dp}.
\end{equation} 

Because $h(p \mid x)$ is generally not available in closed form, the favoured types 
of priors until the introduction of Markov chain Monte Carlo methods have been 
those allowing explicit computations, namely ``conjugate priors.'' These are prior 
distributions for which the corresponding posterior distributions are themselves 
members of the original prior family, the Bayesian updating being accomplished 
through updating of parameters. For a realization $x$ of $X \sim \text{BIN}(n,p)$, 
a family of conjugate priors is the family of beta distributions 
$\text{BETA}(\alpha,\beta)$ where we note from equation (\ref{E:2}) that for $x = 
0,1, \ldots ,n$,
\begin{align*}
  h(p \mid x) &= 
\frac{\binom{n}{x}p^{x}(1-p)^{n-x}\frac{\Gamma(\alpha+\beta)}
{\Gamma(\alpha)\Gamma(\beta)}p^{\alpha-1}(1-p)^{\beta-1}}
{\int_{0}^{1}\binom{n}{x}p^{x}(1-p)^{n-x}\frac{\Gamma(\alpha+\beta)}
{\Gamma(\alpha)\Gamma(\beta)}p^{\alpha-1}(1-p)^{\beta-1} \, dp} \\
              &= 
\frac{\Gamma(\alpha+\beta+n)}{\Gamma(\alpha+x)\Gamma(\beta+n-x)}
p^{\alpha+x-1}(1-p)^{\beta+n-x-1}, \quad 0 \leq p \leq 1, \\
              &= 0, \quad \text{otherwise.}
\end{align*}
That is, the posterior distribution is also beta with updated parameters 
$\alpha+x$ and $\beta+n-x$.

\begin{center}
  \textbf{3. \quad RELATIONSHIPS AND PROPERTIES OF THE POISSON DISTRIBUTION}
\end{center}

Suppose that random variable $X$ has a Poisson distribution with parameter $\lambda 
\geq 0$, denoted by $X \sim \text{POI}(\lambda)$. Then, for $\lambda > 0$, the pmf 
of $X$, denoted by $f_{X}(x \mid \lambda)$, is
\begin{align*}
  f_{X}(x \mid \lambda) = P(X = x \mid \lambda) &= 
\frac{e^{-\lambda}\lambda^{x}}{x!}, \quad x = 0,1,2, \ldots , \\
                                                &= 0, \quad \text{otherwise,}
\end{align*}
and $f_{X}(0 \mid 0) = 1$.

Suppose random variable $Y$ has a gamma distribution with parameters $\alpha > 0$ 
and $\beta > 0$, denoted by $Y \sim \text{GAM}(\alpha,\beta)$. Then the pdf of $Y$, 
denoted \enlargethispage{\baselineskip} by $f_{Y}(y \mid \alpha,\beta)$, is 
\begin{align*}
  f_{Y}(y \mid \alpha,\beta) &= \frac{1}{\beta^{\alpha}\Gamma(\alpha)}
y^{\alpha-1}e^{-y/\beta}, \quad y > 0, \\
                             &= 0, \quad \text{otherwise.}
\end{align*}

Successive integration by parts leads to a relationship between the cdf's of the 
Poisson and gamma distributions. If $X \sim \text{POI}(\lambda)$ and $Y \sim 
\text{GAM}(i+1,2)$ for nonnegative integer $i$, then
\begin{equation} \label{E:3}
  \sum_{x=0}^{i}\frac{e^{-\lambda}\lambda^{x}}{x!} = 1 - 
\frac{1}{2^{i+1}i!}\int_{0}^{2\lambda}t^{i}e^{-t/2} \, dt.
\end{equation}
That is, $F_{X}(i \mid \lambda) = P(X \leq i \mid \lambda) = 1 - P(Y \leq 2\lambda 
\mid i+1,2) = 1 - F_{Y}(2\lambda \mid i+1,2)$. For fixed nonnegative integer $i$, 
it follows from equation (\ref{E:3}) that the function $P(X \leq i \mid \lambda)$ 
is continuous and decreasing in $\lambda$; for positive integer $j$, 
$P(X \geq j \mid \lambda) = 1 - P(X \leq j-1 \mid \lambda)$ is continuous and 
increasing in $\lambda$; and for $1 \leq i \leq j$, $P(i \leq X \leq j \mid \lambda)$ 
is continuous, and increasing for $0 \leq \lambda < \lambda(i,j)$ and decreasing 
for $\lambda \geq \lambda(i,j)$ with maximum at $\lambda = \lambda(i,j) = 
(i \cdots j)^{1/(j-i+1)}$. Also, $\lambda(0,j) = 0$ for $j \geq 0$.

Suppose that the Poisson parameter $\lambda$ is unknown and we wish to estimate it 
using Bayesian methods. For a realization $x$ of $X \sim \text{POI}(\lambda)$, a 
family of conjugate priors is the family of gamma distributions $\text{GAM}(\alpha,
\beta)$ where for $x = 0,1,2, \cdots $, the pdf $h(\lambda \mid x)$ of the posterior 
distribution is given by 
\begin{align*}
  h(\lambda \mid x) &= 
\frac{\frac{e^{-\lambda}\lambda^{x}}{x!}\frac{1}{\beta^{\alpha}\Gamma(\alpha)}
\lambda^{\alpha-1}e^{-\lambda/\beta}}{\int_{0}^{\infty}\frac{e^{-\lambda}\lambda^{x}}
{x!}\frac{1}{\beta^{\alpha}\Gamma(\alpha)}\lambda^{\alpha-1}e^{-\lambda/\beta} \, 
d\lambda} \\
                    &= \frac{1}{[\beta/(1+\beta)]^{\alpha+x}\Gamma(\alpha+x)}
\lambda^{\alpha+x-1}e^{-\lambda/[\beta/(1+\beta)]}, \quad \lambda > 0, \\
                    &= 0, \quad \text{otherwise.}
\end{align*}
That is, the posterior distribution is also gamma with updated parameters $\alpha + 
x$ and $\beta/(1 + \beta)$.

\pagebreak

\begin{center}
  \textbf{4. \quad RELATIONSHIPS AND PROPERTIES OF THE HYPERGEOMETRIC DISTRIBUTION}
\end{center}

Suppose that integer-valued random variable $X$ has a hypergeometric distribution 
with parameters $n$, $M$, and $N$, denoted by $X \sim \text{HYP}(n,M,N)$, where $n$, 
$M$, and $N$  are integers with $1 \leq n \leq N$ and $0 \leq M \leq N$. Then, for 
given $n$ and $N$, and for $0 < M < N$, the pmf of $X$, denoted by 
$f_{X}(x \mid M)$, is
\begin{align} \label{E:4}
  f_{X}(x \mid M) = P(X = x \mid M) &= \frac{\binom{M}{x}\binom{N-M}{n-x}}
{\binom{N}{n}}, \quad \text{max}(0,n-N+M) \leq x \leq \text{min}(n,M), \notag \\
                                    &= 0, \quad \text{otherwise,}
\end{align}
and $f_{X}(0 \mid 0) = f_{X}(n \mid N) = 1$.  

Suppose that random variable $Y$ has a specially defined discrete distribution with 
parameters $a$, $b$, and $c$, denoted by $Y \sim \text{ABC}(a,b,c)$, where $a$, $b$, 
and $c$ are nonnegative integers. Then, for $c > 0$, the pmf of $Y$, denoted by 
$f_{Y}(y \mid a,b,c)$, is
\begin{align*}
  f_{Y}(y \mid a,b,c) = P(Y = y \mid a,b,c) &= \frac{\binom{a+y}{a}\binom{b+c-y}{b}}
{\binom{a+b+c+1}{a+b+1}}, \quad y = 0,1, \ldots ,c, \\
                                            &= 0, \quad \text{otherwise,}
\end{align*}
and $f_{Y}(0 \mid a,b,0) = 1$. We note that formula (12.16) of Feller (1968, p.65) 
can be used to prove that

\smallskip 
\begin{equation*}
  \sum_{y=0}^{c}\binom{a+y}{a}\binom{b+c-y}{b} = \binom{a+b+c+1}{a+b+1}.
\end{equation*}
\smallskip

We also note that the ABC distribution is just a special case of the Polya 
(or beta-binomial) distribution (Dyer and Pierce, 1993, p.2130). From equation 
(\ref{E:4}), it easily follows that $P(X \leq n \mid M) = 1$ for $0 \leq M \leq N$.											
For $0 \leq i < n \leq N$ and $0 \leq M \leq N$, we have from equation (\ref{E:4}) 
that 

\begin{align*}
  \binom{N}{n}P(X \leq i \mid M) &= \sum_{x=0}^{i}\binom{M}{x}\binom{N-M}{n-x} 
\notag \\
                                 &= \sum_{x=0}^{i}\binom{M}{x}
\left[ \binom{N-M-1}{n-x-1} + \binom{N-M-1}{n-x} \right] \notag \\
                                 &= \sum_{x=0}^{i}\binom{M}{x}\binom{N-M-1}{n-x-1} 
+ \sum_{x=0}^{i}\binom{M}{x}\binom{N-M-1}{n-x} \notag \\
                                 &= \sum_{x=1}^{i+1}\binom{M}{x-1}\binom{N-M-1}{n-x} 
+ \sum_{x=0}^{i}\binom{M}{x}\binom{N-M-1}{n-x} \notag \\
\end{align*}

\pagebreak

\begin{align} \label{E:5}
\phantom{\binom{N}{n}P(X \leq i \mid M)} &= \binom{M}{i}\binom{N-M-1}{n-i-1} - 
\binom{M}{-1}\binom{N-M-1}{n} \notag \\
                                         & \qquad  + \sum_{x=0}^{i}\left[ \binom{M}{x-1} + 
\binom{M}{x} \right]\binom{N-M-1}{n-x} \notag \\
                                         &= \binom{M}{i}\binom{N-M-1}{n-i-1} + 
\sum_{x=0}^{i}\binom{M+1}{x}\binom{N-M-1}{n-x} \notag \\
                                         &= \binom{M}{i}\binom{N-M-1}{n-i-1} + 
\binom{N}{n}P(X \leq i \mid M+1), 
\end{align}
where by definition $\binom{M}{-1} = 0$, $\binom{M}{i} = 0$ if $M < i$, and 
$\binom{N-M-1}{n-i-1} = 0$ if $M > N-n+i$. Furthermore, from the recursion relationship 
in equation (\ref{E:5}), it follows that 
\begin{align} \label{E:6}
  P(X \leq i \mid M) &= \sum_{k=M}^{N-n+i}\binom{k}{i}\binom{N-k-1}{n-i-1}\biggr/
\binom{N}{n} \notag \\
                     &= \sum_{k=M-i}^{N-n}\binom{i+k}{i}\binom{n-i-1+N-n-k}{n-i-1}
\biggr/\binom{N}{n} \notag \\
                     &= 1 - \sum_{k=0}^{M-i-1}\binom{i+k}{i}\binom{n-i-1+N-n-k}{n-i-1}
\biggr/\binom{N}{n}.
\end{align}
That is, if $X \sim \text{HYP}(n,M,N)$ and $Y \sim \text{ABC}(i,n-i-1,N-n)$ for 
integer $i$, $0 \leq i < n \leq N$, then $F_{X}(i \mid M) = P(X \leq i \mid M) = 1 
- P(Y \leq M-i-1 \mid i,n-i-1,N-n) = 1 - F_{Y}(M-i-1 \mid i,n-i-1,N-n)$ where, in 
particular, 
\begin{align} \label{E:7}
  P(X \leq i \mid M) &= 1, \quad \text{if} \quad 0 \leq M \leq i, \notag \\
                     &= 0, \quad \text{if} \quad N-n+i < M \leq N.
\end{align}

For $0 < i \leq j < n \leq N$ and $0 \leq M \leq N$, we have from equation 
(\ref{E:5}) that 
\begin{align} \label{E:8}
  \binom{N}{n}P(i \leq X \leq j \mid M) &= \binom{N}{n}P(X \leq j \mid M) - 
\binom{N}{n}P(X \leq i-1 \mid M) \notag \\
                                        &= \binom{M}{j}\binom{N-M-1}{n-j-1} + 
\binom{N}{n}P(X \leq j \mid M+1) \notag \\
                                        & \qquad - \binom{M}{i-1}\binom{N-M-1}{n-i} 
- \binom{N}{n}P(X \leq i-1 \mid M+1) \notag \\
                                        &= \binom{M}{j}\binom{N-M-1}{n-j-1} - 
\binom{M}{i-1}\binom{N-M-1}{n-i} \notag \\
                                        & \qquad + \binom{N}{n}P(i \leq X \leq j 
\mid M+1).
\end{align}

Similar to the determination of equation (\ref{E:6}), it follows from the recursion 
relationship in equation (\ref{E:8}) that
\begin{align} \label{E:9}
  P(i \leq X \leq j \mid M) &= \sum_{k=M}^{N-n+j}\binom{k}{j}\binom{N-k-1}{n-j-1}
\biggr/\binom{N}{n} - \sum_{l=M}^{N-n+i-1}\binom{l}{i-1}\binom{N-l-1}{n-i}
\biggr/\binom{N}{n} \notag \\
                            &= \sum_{k=M-j}^{N-n}\binom{j+k}{j}\binom{n-j-1+N-n-k}
{n-j-1}\biggr/\binom{N}{n} \notag \\
                            & \qquad - \sum_{l=M-i+1}^{N-n}\binom{i-1+l}{i-1}
\binom{n-i+N-n-l}{n-i}\biggr/\binom{N}{n} \notag \\
                            &= \sum_{l=0}^{M-i}\binom{i-1+l}{i-1}\binom{n-i+N-n-l}
{n-i}\biggr/\binom{N}{n} \notag \\
                            & \qquad - \sum_{k=0}^{M-j-1}\binom{j+k}{j}
\binom{n-j-1+N-n-k}{n-j-1}\biggr/\binom{N}{n}
\end{align}
where, in particular,
\begin{equation} \label{E:10}
  P(i \leq X \leq j \mid M) = 0, \quad \text{if either} \quad 0 \leq M < i \quad 
\text{or} \quad N-n+j < M \leq N.
\end{equation} 
\noindent
We note in equation (\ref{E:8}) that the difference
\begin{align} \label{E:11}
  \binom{M}{j}\binom{N-M-1}{n-j-1} - \binom{M}{i-1}\binom{N-M-1}{n-i} &=
-\binom{N-i}{n-i} < 0, \quad \text{if} \quad M=i-1, \notag \\
                                                                      &=
\binom{N-n+j}{j} > 0, \quad \text{if} \quad M=N-n+j,
\end{align}
and for $i \leq M < N-n+j$, the same difference 
\begin{align} \label{E:12}
& \binom{M}{j}\binom{N-M-1}{n-j-1} - \binom{M}{i-1}\binom{N-M-1}{n-i} \notag \\ 
&= \frac{M!}{j!(M-j)!}\frac{(N-M-1)!}{(n-j-1)!(N-M-n+j)!}  
- \frac{M!}{(i-1)!(M-i+1)!}\frac{(N-M-1)!}{(n-i)!(N-M-n+i-1)!} \notag \\                                                                      
&= \frac{M!(N-M-1)!}{(i-1)!(M-j)!(n-j-1)!(N-M-n+i-1)!} \notag \\
& \qquad \times \left[ \frac{1}{(j \cdots i)}\frac{1}{(N-M-n+j) \cdots (N-M-n+i)} 
\right. \notag \\                                                                        
& \qquad \qquad - \left. \frac{1}{(M-i+1) \cdots (M-j+1)}\frac{1}{(n-i) \cdots 
(n-j)} \right]
\end{align}
where as $M$ increases, the term $1/(N-M-n+j) \cdots (N-M-n+i)$ increases and the 
term $1/(M-i+1) \cdots (M-j+1)$ decreases so that as $M$ increases between $i-1$ 
and $N-n+j$, the difference $\binom{M}{j}\binom{N-M-1}{n-j-1} - \binom{M}{i-1}
\binom{N-M-1}{n-i}$ goes from being negative to being positive and staying positive.
											
In summary, $P(X \leq n \mid M)$ equals 1 for $0 \leq M \leq N$, and for fixed 
integer $i$, $0 \leq i < n \leq N$, we see from equations (\ref{E:6}) and (\ref{E:7}) 
that $P(X \leq i \mid M )$ equals 1 for $0 \leq M \leq i$, is decreasing for 
$i < M \leq N-n+i$, and equals 0 for $N-n+i < M \leq N$; $P(X \geq n+1 \mid M)$ 
equals 0 for $0 \leq M \leq N$, and for fixed integer $j$, $1 \leq j \leq n \leq N$, 
$P(X \geq j \mid M) = 1 - P(X \leq j-1 \mid M)$ equals 0 for $0 \leq M \leq j-1$, is 
increasing for $j-1 < M \leq N-n+j-1$, and equals 1 for $N-n+j-1 < M \leq N$; and 
we see from equations (\ref{E:8}) to (\ref{E:12}) that for fixed integers $i$ and $j$, 
$0 < i \leq j < n \leq N$ where we define
\begin{equation*}
  M_{n,N}(i,j) = \text{min}\{ M \mid i \leq M \leq N-n+j \quad \text{and} 
\quad 
\textstyle{\binom{M}{j}\binom{N-M-1}{n-j-1} \geq  \binom{M}{i-1}\binom{N-M-1}{n-i}} \},
\end{equation*}
$P(i \leq X \leq j \mid M)$ equals 0 for $0 \leq M < i$, is increasing for 
$i \leq M < M_{n,N}(i,j)$, is decreasing for $M_{n,N}(i,j) + 1 < M \leq N-n+j$, and 
equals 0 for $N-n+j < M \leq N$ with maximum at either $M_{n,N}(i,j)$ if   
$\binom{M}{j}\binom{N-M-1}{n-j-1} > \binom{M}{i-1}\binom{N-M-1}{n-i}$ for 
$M = M_{n,N}(i,j)$ so that $P(i \leq X \leq j \mid M_{n,N}(i,j)) > 
P(i \leq X \leq j \mid M_{n,N}(i,j) + 1)$ or maximum at both $M_{n,N}(i,j)$ and 
$M_{n,N}(i,j) + 1$ if $\binom{M}{j}\binom{N-M-1}{n-j-1} = \binom{M}{i-1}\binom{N-M-1}{n-i}$
for $M = M_{n,N}(i,j)$ so that  $P(i \leq X \leq j \mid M_{n,N}(i,j)) = 
P(i \leq X \leq j \mid M_{n,N}(i,j) + 1)$.

Suppose that the hypergeometric parameters $n$ and $N$ are known but $M$ is not and 
we wish to estimate it using Bayesian methods. For a realization $x$ of $X \sim 
\text{HYP}(n,M,N)$, a family of conjugate priors for $M - x$ is the family of discrete 
distributions $\text{ABC}(a,b,N)$ where for $x = 0,1, \ldots ,n$, the pmf $h(M \mid x)$ of 
the posterior distribution for $M$ is given by
\begin{align} \label{E:13}
  h(M \mid x) &= \frac{\frac{\binom{M}{x}\binom{N-M}{n-x}}{\binom{N}{n}}
\frac{\binom{a+M}{a}\binom{b+N-M}{b}}{\binom{a+b+N+1}{a+b+1}}}
{\sum_{M=x}^{N-n+x}\frac{\binom{M}{x}\binom{N-M}{n-x}}{\binom{N}{n}}
\frac{\binom{a+M}{a}\binom{b+N-M}{b}}{\binom{a+b+N+1}{a+b+1}}} \notag \\
              &= \frac{\binom{a+M}{a+x}\binom{b+N-M}{b+n-x}}
{\binom{a+b+N+1}{a+b+n+1}}, \quad x \leq M \leq N-n+x, \notag \\
              &= 0, \quad \text{otherwise,}
\end{align}                              											
from which it easily follows that the pmf $h(M-x \mid x)$ of the posterior distribution 
for $M-x$ is given by 
\begin{align} \label{E:14}
  h(M-x \mid x) &= \frac{\binom{a+x+M-x}{a+x}\binom{b+n-x+N-n-M+x}{b+n-x}}
{\binom{a+x+b+n-x+N-n+1}{a+x+b+n-x+1}}, \quad 0 \leq M-x \leq N-n, \notag \\
                &= 0, \quad \text{otherwise.}
\end{align}
That is, the posterior distribution for $M-x$ is also ABC with updated parameters 
$a+x$, $b+n-x$, and $N-n$.

Finally, we note that as a family of conjugate priors for the hypergeometric 
distribution 
\linebreak 
$\text{HYP}(n,M,N)$, the family of discrete distributions $\text{ABC}(a,b,N)$ has, 
in addition to unimodal members, strictly increasing members $\text{ABC}(a,0,N)$, strictly 
decreasing members $\text{ABC}(0,b,N)$, and the discrete uniform distribution 
$\text{ABC}(0,0,N)$.

\newpage
\begin{center}
  \textbf{REFERENCES}
\end{center}
\noindent
Dyer, D. and Pierce, R. L. (1993), "On the choice of the prior distribution in hypergeometric \\
\indent sampling," \emph{Communications in Statistics - Theory and Methods}, 22(8), 2125-2146. \\
\noindent
Feller, W. (1968), \emph{An Introduction to Probability Theory and Its Applications}, 
Vol.1, (3rd ed.), \\
\indent John Wiley \& Sons, Inc.
											
\end{document}